\documentclass[reqno,12pt]{amsart}
\usepackage[latin1]{inputenc}
\usepackage[T1]{fontenc}
\usepackage{texdraw}
\newtheorem{theorem}{Theorem}%[section]
\newtheorem{lemma}[theorem]{Lemma}

\theoremstyle{definition}

\theoremstyle{remark}

\newtheorem*{remark}{Remark}

\def\({\left(}
\def\){\right)}
\def\[{\left[}
\def\]{\right]}

\begin{document}
\thispagestyle{empty}

\catcode`\@=11
\font\tenln    = line10
\font\tenlnw   = linew10

\thinlines
\newskip\Einheit \Einheit=0.6cm
\newcount\xcoord \newcount\ycoord
\newdimen\xdim \newdimen\ydim \newdimen\PfadD@cke \newdimen\Pfadd@cke
\PfadD@cke1pt \Pfadd@cke0.5pt
\def\PfadDicke#1{\PfadD@cke#1 \divide\PfadD@cke by2 \Pfadd@cke\PfadD@cke \multiply\PfadD@cke by2}
\long\def\LOOP#1\REPEAT{\def\BODY{#1}\ITERATE}
\def\ITERATE{\BODY \let\next\ITERATE \else\let\next\relax\fi \next}
\let\REPEAT=\fi
\def\Punkt{\hbox{\raise-2pt\hbox to0pt{\hss\scriptsize$\bullet$\hss}}}
\def\DuennPunkt(#1,#2){\unskip
  \raise#2 \Einheit\hbox to0pt{\hskip#1 \Einheit
          \raise-2.5pt\hbox to0pt{\hss\normalsize$\bullet$\hss}\hss}}
\def\NormalPunkt(#1,#2){\unskip
  \raise#2 \Einheit\hbox to0pt{\hskip#1 \Einheit
          \raise-3pt\hbox to0pt{\hss\large$\bullet$\hss}\hss}}
\def\DickPunkt(#1,#2){\unskip
  \raise#2 \Einheit\hbox to0pt{\hskip#1 \Einheit
          \raise-4pt\hbox to0pt{\hss\Large$\bullet$\hss}\hss}}
\def\Kreis(#1,#2){\unskip
  \raise#2 \Einheit\hbox to0pt{\hskip#1 \Einheit
          \raise-4pt\hbox to0pt{\hss\Large$\circ$\hss}\hss}}
\def\Diagonale(#1,#2)#3{\unskip\leavevmode
  \xcoord#1\relax \ycoord#2\relax
      \raise\ycoord \Einheit\hbox to0pt{\hskip\xcoord \Einheit
         \unitlength\Einheit
         \line(1,1){#3}\hss}}
\def\AntiDiagonale(#1,#2)#3{\unskip\leavevmode
  \xcoord#1\relax \ycoord#2\relax \advance\xcoord by -0.05\relax
      \raise\ycoord \Einheit\hbox to0pt{\hskip\xcoord \Einheit
         \unitlength\Einheit
         \line(1,-1){#3}\hss}}
\def\Pfad(#1,#2),#3\endPfad{\unskip\leavevmode
  \xcoord#1 \ycoord#2 \thicklines\ZeichnePfad#3\endPfad\thinlines}
\def\ZeichnePfad#1{\ifx#1\endPfad\let\next\relax
  \else\let\next\ZeichnePfad
    \ifnum#1=1
      \raise\ycoord \Einheit\hbox to0pt{\hskip\xcoord \Einheit
         \vrule height\Pfadd@cke width1 \Einheit depth\Pfadd@cke\hss}%
      \advance\xcoord by 1
    \else\ifnum#1=2
      \raise\ycoord \Einheit\hbox to0pt{\hskip\xcoord \Einheit
        \hbox{\hskip-1pt\vrule height1 \Einheit width\PfadD@cke depth0pt}\hss}%
      \advance\ycoord by 1
    \else\ifnum#1=3
      \raise\ycoord \Einheit\hbox to0pt{\hskip\xcoord \Einheit
         \unitlength\Einheit
         \line(1,1){1}\hss}
      \advance\xcoord by 1
      \advance\ycoord by 1
    \else\ifnum#1=4
      \raise\ycoord \Einheit\hbox to0pt{\hskip\xcoord \Einheit
         \unitlength\Einheit
         \line(1,-1){1}\hss}
      \advance\xcoord by 1
      \advance\ycoord by -1
    \else\ifnum#1=5  
\raise\ycoord \Einheit\hbox to0pt{\hskip\xcoord \Einheit
         \unitlength\Einheit
         \line(0,-1){1}\hss}
      \advance\ycoord by -1
    \fi\fi\fi\fi\fi
  \fi\next}
\def\hSSchritt{\leavevmode\raise-.4pt\hbox to0pt{\hss.\hss}\hskip.2\Einheit
  \raise-.4pt\hbox to0pt{\hss.\hss}\hskip.2\Einheit
  \raise-.4pt\hbox to0pt{\hss.\hss}\hskip.2\Einheit
  \raise-.4pt\hbox to0pt{\hss.\hss}\hskip.2\Einheit
  \raise-.4pt\hbox to0pt{\hss.\hss}\hskip.2\Einheit}
\def\vSSchritt{\vbox{\baselineskip.2\Einheit\lineskiplimit0pt
\hbox{.}\hbox{.}\hbox{.}\hbox{.}\hbox{.}}}
\def\DSSchritt{\leavevmode\raise-.4pt\hbox to0pt{%
  \hbox to0pt{\hss.\hss}\hskip.2\Einheit
  \raise.2\Einheit\hbox to0pt{\hss.\hss}\hskip.2\Einheit
  \raise.4\Einheit\hbox to0pt{\hss.\hss}\hskip.2\Einheit
  \raise.6\Einheit\hbox to0pt{\hss.\hss}\hskip.2\Einheit
  \raise.8\Einheit\hbox to0pt{\hss.\hss}\hss}}
\def\dSSchritt{\leavevmode\raise-.4pt\hbox to0pt{%
  \hbox to0pt{\hss.\hss}\hskip.2\Einheit
  \raise-.2\Einheit\hbox to0pt{\hss.\hss}\hskip.2\Einheit
  \raise-.4\Einheit\hbox to0pt{\hss.\hss}\hskip.2\Einheit
  \raise-.6\Einheit\hbox to0pt{\hss.\hss}\hskip.2\Einheit
  \raise-.8\Einheit\hbox to0pt{\hss.\hss}\hss}}
\def\SPfad(#1,#2),#3\endSPfad{\unskip\leavevmode
  \xcoord#1 \ycoord#2 \ZeichneSPfad#3\endSPfad}
\def\ZeichneSPfad#1{\ifx#1\endSPfad\let\next\relax
  \else\let\next\ZeichneSPfad
    \ifnum#1=1
      \raise\ycoord \Einheit\hbox to0pt{\hskip\xcoord \Einheit
         \hSSchritt\hss}%
      \advance\xcoord by 1
    \else\ifnum#1=2
      \raise\ycoord \Einheit\hbox to0pt{\hskip\xcoord \Einheit
        \hbox{\hskip-2pt \vSSchritt}\hss}%
      \advance\ycoord by 1
    \else\ifnum#1=3
      \raise\ycoord \Einheit\hbox to0pt{\hskip\xcoord \Einheit
         \DSSchritt\hss}
      \advance\xcoord by 1
      \advance\ycoord by 1
    \else\ifnum#1=4
      \raise\ycoord \Einheit\hbox to0pt{\hskip\xcoord \Einheit
         \dSSchritt\hss}
      \advance\xcoord by 1
      \advance\ycoord by -1
    \fi\fi\fi\fi
  \fi\next}
\def\Koordinatenachsen(#1,#2){\unskip
 \hbox to0pt{\hskip-.5pt\vrule height#2 \Einheit width.5pt depth1 \Einheit}%
 \hbox to0pt{\hskip-1 \Einheit \xcoord#1 \advance\xcoord by1
    \vrule height0.25pt width\xcoord \Einheit depth0.25pt\hss}}
\def\Koordinatenachsen(#1,#2)(#3,#4){\unskip
 \hbox to0pt{\hskip-.5pt \ycoord-#4 \advance\ycoord by1
    \vrule height#2 \Einheit width.5pt depth\ycoord \Einheit}%
 \hbox to0pt{\hskip-1 \Einheit \hskip#3\Einheit 
    \xcoord#1 \advance\xcoord by1 \advance\xcoord by-#3 
    \vrule height0.25pt width\xcoord \Einheit depth0.25pt\hss}}
\def\Gitter(#1,#2){\unskip \xcoord0 \ycoord0 \leavevmode
  \LOOP\ifnum\ycoord<#2
    \loop\ifnum\xcoord<#1
      \raise\ycoord \Einheit\hbox to0pt{\hskip\xcoord \Einheit\Punkt\hss}%
      \advance\xcoord by1
    \repeat
    \xcoord0
    \advance\ycoord by1
  \REPEAT}
\def\Gitter(#1,#2)(#3,#4){\unskip \xcoord#3 \ycoord#4 \leavevmode
  \LOOP\ifnum\ycoord<#2
    \loop\ifnum\xcoord<#1
      \raise\ycoord \Einheit\hbox to0pt{\hskip\xcoord \Einheit\Punkt\hss}%
      \advance\xcoord by1
    \repeat
    \xcoord#3
    \advance\ycoord by1
  \REPEAT}
\def\Label#1#2(#3,#4){\unskip \xdim#3 \Einheit \ydim#4 \Einheit
  \def\lo{\advance\xdim by-.5 \Einheit \advance\ydim by.5 \Einheit}%
  \def\llo{\advance\xdim by-.25cm \advance\ydim by.5 \Einheit}%
  \def\loo{\advance\xdim by-.5 \Einheit \advance\ydim by.25cm}%
  \def\o{\advance\ydim by.25cm}%
  \def\ro{\advance\xdim by.5 \Einheit \advance\ydim by.5 \Einheit}%
  \def\rro{\advance\xdim by.25cm \advance\ydim by.5 \Einheit}%
  \def\roo{\advance\xdim by.5 \Einheit \advance\ydim by.25cm}%
  \def\l{\advance\xdim by-.30cm}%
  \def\r{\advance\xdim by.30cm}%
  \def\lu{\advance\xdim by-.5 \Einheit \advance\ydim by-.6 \Einheit}%
  \def\llu{\advance\xdim by-.25cm \advance\ydim by-.6 \Einheit}%
  \def\luu{\advance\xdim by-.5 \Einheit \advance\ydim by-.30cm}%
  \def\u{\advance\ydim by-.30cm}%
  \def\ru{\advance\xdim by.5 \Einheit \advance\ydim by-.6 \Einheit}%
  \def\rru{\advance\xdim by.25cm \advance\ydim by-.6 \Einheit}%
  \def\ruu{\advance\xdim by.5 \Einheit \advance\ydim by-.30cm}%
  #1\raise\ydim\hbox to0pt{\hskip\xdim
     \vbox to0pt{\vss\hbox to0pt{\hss$#2$\hss}\vss}\hss}%
}
\catcode`\@=12

%Christian's texdraw definitions
\def\ldreieck{\bsegment
  \rlvec(0.866025403784439 .5) \rlvec(0 -1)
  \rlvec(-0.866025403784439 .5)  
  \savepos(0.866025403784439 -.5)(*ex *ey)
        \esegment
  \move(*ex *ey)
        }
\def\rdreieck{\bsegment
  \rlvec(0.866025403784439 -.5) \rlvec(-0.866025403784439 -.5)  \rlvec(0 1)
  \savepos(0 -1)(*ex *ey)
        \esegment
  \move(*ex *ey)
        }
\def\rhombus{\bsegment
  \rlvec(0.866025403784439 .5) \rlvec(0.866025403784439 -.5) 
  \rlvec(-0.866025403784439 -.5)  \rlvec(0 1)        
  \rmove(0 -1)  \rlvec(-0.866025403784439 .5) 
  \savepos(0.866025403784439 -.5)(*ex *ey)
        \esegment
  \move(*ex *ey)
        }
\def\RhombusA{\bsegment
  \rlvec(0.866025403784439 .5) \rlvec(0.866025403784439 -.5) 
  \rlvec(-0.866025403784439 -.5) \rlvec(-0.866025403784439 .5) 
  \savepos(0.866025403784439 -.5)(*ex *ey)
        \esegment
  \move(*ex *ey)
        }
\def\RhombusB{\bsegment
  \rlvec(0.866025403784439 .5) \rlvec(0 -1)
  \rlvec(-0.866025403784439 -.5) \rlvec(0 1) 
  \savepos(0 -1)(*ex *ey)
        \esegment
  \move(*ex *ey)
        }
\def\RhombusC{\bsegment
  \rlvec(0.866025403784439 -.5) \rlvec(0 -1)
  \rlvec(-0.866025403784439 .5) \rlvec(0 1) 
  \savepos(0.866025403784439 -.5)(*ex *ey)
        \esegment
  \move(*ex *ey)
        }
\def\hdSchritt{\bsegment 
  \lpatt(.05 .13)
  \rlvec(0.866025403784439 -.5) 
  \savepos(0.866025403784439 -.5)(*ex *ey)
        \esegment
  \move(*ex *ey)
        }
\def\vdSchritt{\bsegment
  \lpatt(.05 .13)
  \rlvec(0 -1) 
  \savepos(0 -1)(*ex *ey)
        \esegment
  \move(*ex *ey)
        }

\def\ringerl(#1 #2){\move(#1 #2)\fcir f:0 r:.15}
\def\knoten{\bsegment \fcir f:0 r:.15 \esegment}

\def\hantel(#1 #2){\fcir f:0 r:.1 \rlvec(#1 #2) \fcir f:0 r:.1}

\def\hex{\bsegment
        \rlvec(1 0)  \rlvec(.5 -.866025403784439) \rlvec(-.5 -.866025403784439)
        \rlvec(-1 0) \rlvec(-.5 .866025403784439) \rlvec(.5 .866025403784439)
        \savepos(1.5 -.866025403784439)(*ex *ey)
         \esegment
        \move(*ex *ey)
}
\def\hexa{\bsegment \lcir r:.1
        \rlvec(1 0) \fcir f:0 r:.2  \rlvec(.5 -.866025403784439)  
         \lcir r:.1 \rlvec(-.5 -.866025403784439) \fcir f:0 r:.2 
        \rlvec(-1 0) 
       \lcir r:.1 \rlvec(-.5 .866025403784439) 
        \fcir f:0 r:.2 \rlvec(.5 .866025403784439)
        \savepos(1.5 -.866025403784439)(*ex *ey)
         \esegment
        \move(*ex *ey)
}
\def\hexb{\bsegment \fcir f:0 r:.2
        \ravec(1 0)   \lcir r:.1 \rmove(.5 -.866025403784439) 
\ravec(-.5 .866025403784439) \rmove(.5 -.866025403784439) 
         \fcir f:0 r:.2 \ravec(-.5 -.866025403784439) \lcir r:.1 
        \rmove(-1 0)\ravec(1 0)\rmove(-1 0) 
       \fcir f:0 r:.2 \ravec(-.5 .866025403784439) 
        \lcir r:.1 \rmove(.5 .866025403784439)
\ravec(-.5 -.866025403784439)\rmove(.5 .866025403784439)
        \savepos(1.5 -.866025403784439)(*ex *ey)
         \esegment
        \move(*ex *ey)
}

\def\shex{\bsegment 
        \lpatt(.05 .13)
        \rlvec(1 0)  \rlvec(.5 -.866025403784439) \rlvec(-.5 -.866025403784439)
        \rlvec(-1 0) \rlvec(-.5 .866025403784439) \rlvec(.5 .866025403784439)
        \savepos(1.5 -.866025403784439)(*ex *ey)
         \esegment
        \move(*ex *ey)
}

\def\RhombiA{\bsegment
  \rlvec(0.866025403784439 .5) \rlvec(0.866025403784439 -.5) 
  
  \rlvec(-0.866025403784439 -.5) \rlvec(-0.866025403784439 .5)
  \lfill f:0.7 
  \savepos(0.866025403784439 -.5)(*ex *ey)
        \esegment
  \move(*ex *ey)
        }
\def\RhombiB{\bsegment
  \rlvec(0.866025403784439 .5) \rlvec(0 -1)
  \rlvec(-0.866025403784439 -.5) \rlvec(0 1) \lfill f:.9
  \savepos(0 -1)(*ex *ey)
        \esegment
  \move(*ex *ey)
        }
\def\RhombiC{\bsegment
  \rlvec(0.866025403784439 -.5) \rlvec(0 -1)
  \rlvec(-0.866025403784439 .5) \rlvec(0 1) \lfill f:.3
  \savepos(0.866025403784439 -.5)(*ex *ey)
        \esegment
  \move(*ex *ey)
        }

\def\RhA{\bsegment
  \rlvec(0.866025403784439 .5) \rlvec(0.866025403784439 -.5) 
  \lfill f:0.5
  \rlvec(-0.866025403784439 -.5) \rlvec(-0.866025403784439 .5)
  \lfill f:0.5 
  \savepos(1.732 0)(*ex *ey)
        \esegment
  \move(*ex *ey)
        }
\def\RhB{\bsegment
\rlvec(0 1) \rlvec(0.866025403784439 .5) \lfill f:.5
\rlvec(0 -1) \rlvec(-0.866025403784439 -.5) \lfill f:.5
  \savepos(0.866025 1.5)(*ex *ey)
        \esegment
  \move(*ex *ey)
        }
\def\RhC{\bsegment
  \rlvec(0.866025403784439 -.5) \rlvec(0 -1) \lfill f:.5
  \rlvec(-0.866025403784439 .5) \rlvec(0 1) \lfill f:.5
  \savepos(0.866025403784439 -1.5)(*ex *ey)
        \esegment
  \move(*ex *ey)
        }

\def\nwdom{\bsegment
\rlvec(1 1) \rlvec(2 -2)\rlvec(-1 -1) \rlvec(-2 2)
\savepos(2 0)(*ex *ey)
\esegment
\move(*ex *ey)
} 

\def\nedom{\bsegment
\rlvec(1 1) \rlvec(1 -1)\rlvec(-2 -2) \rlvec(-1 1) \rlvec(1 1)
\savepos(2 0)(*ex *ey)
\esegment
\move(*ex *ey)
}

\def\fsq{\bsegment \lpatt(.05 .13) \rlvec(1 1)
\rlvec(1 -1)
\rlvec(-1 -1)
\rlvec(-1 1)\lfill f:.8 \esegment}

\def\quadi{\bsegment \rlvec(1 1) \rlvec(1 -1) 
\rlvec(-1 -1) \rlvec(-1 1)
\savepos(2 0)(*ex *ey)
\esegment
\move(*ex *ey)}

        \def\squadi{\bsegment
        \lpatt(.05 .13)
        \rlvec(1 1) \rlvec(1 -1) \rlvec(-1 -1)  \rlvec(-1 1)
        \savepos(2 0)(*ex *ey)
        \esegment
        \move(*ex *ey) }

\def\quadk{\bsegment \knoten \rlvec(1 1) \knoten \rlvec(1 -1) \knoten 
\rlvec(-1 -1) \knoten \rlvec(-1 1)
\savepos(2 0)(*ex *ey)
\esegment
\move(*ex *ey)}

        \def\squadk{\bsegment
        \lpatt(.05 .13) \knoten
        \rlvec(1 1) \knoten \rlvec(1 -1) \knoten \rlvec(-1 -1) \knoten
 \rlvec(-1 1)
        \savepos(2 0)(*ex *ey)
        \esegment
        \move(*ex *ey) }

\def\ringquad{\bsegment
\knoten
\rlvec(1 0) \knoten \rlvec(1 1) \knoten \rlvec(0 1) \knoten 
\rmove(0 -1) 
\rlvec(1 -1) \knoten \rlvec(1 0) \knoten
\rmove(-1 0) 
\rlvec(-1 -1) \knoten \rlvec(0 -1) \knoten
\rmove(0 1) 
\rlvec(-1 1) 
\rmove(1 0)       \lpatt(.05 .13) \lcir r:2 \savepos(4 0)(*ex *ey) \esegment \move(*ex *ey)}
  
\def\urba{\bsegment
\knoten
\rlvec(1 0) \knoten \rlvec(1 1) \knoten \rlvec(0 1) \knoten 
\rmove(0 -1) 
\rlvec(1 -1) \knoten \rlvec(1 0) \knoten
\rmove(-1 0) 
\rlvec(-1 -1) \knoten \rlvec(0 -1) \knoten
\rmove(0 1) 
\rlvec(-1 1) 
\savepos(4 0)(*ex *ey) \esegment \move(*ex *ey)}
  
\def\urbb{\bsegment 
\knoten
\rlvec(1 0) \knoten \lpatt(.05 .13) \rlvec(1 1) \lpatt(1 0) 
\knoten \rlvec(0 1) \knoten 
\rmove(0 -1) 
\lpatt(.05 .13) \rlvec(1 -1) \lpatt(1 0) \knoten \rlvec(1 0) \knoten
\rmove(-1 0) 
\lpatt(.05 .13) \rlvec(-1 -1) \lpatt(1 0) \knoten \rlvec(0 -1) \knoten
\rmove(0 1) 
\lpatt(.05 .13) \rlvec(-1 1) \lpatt(1 0)
\rmove(1 0)  \savepos(4 0)(*ex *ey) \esegment \move(*ex *ey)}

\def\urbc{\bsegment 
\knoten
\rlvec(1 0) \knoten \lpatt(.05 .13) \rlvec(1 1) \lpatt(1 0) 
\knoten \rlvec(0 1) \knoten 
\rmove(0 -1) 
\rlvec(1 -1) \knoten \rlvec(1 0) \knoten
\rmove(-1 0) 
\lpatt(.05 .13) \rlvec(-1 -1) \lpatt(1 0) \knoten \rlvec(0 -1) \knoten
\rmove(0 1) 
\rlvec(-1 1) 
\rmove(1 0)  \savepos(4 0)(*ex *ey) \esegment \move(*ex *ey)}

\def\urbd{\bsegment 
\knoten
\rlvec(1 0) \knoten  \rlvec(1 1)  
\knoten \rlvec(0 1) \knoten 
\rmove(0 -1) 
\lpatt(.05 .13) \rlvec(1 -1) \lpatt(1 0) \knoten \rlvec(1 0) \knoten
\rmove(-1 0) 
\rlvec(-1 -1) \knoten \rlvec(0 -1) \knoten
\rmove(0 1) 
\lpatt(.05 .13) \rlvec(-1 1) \lpatt(1 0)
\rmove(1 0)  \savepos(4 0)(*ex *ey) \esegment \move(*ex *ey)}

\def\dblne#1{\bsegment
\rmove(-.05 .05) \rlvec(#1 #1) \rmove(.1 -.1) \rlvec(-#1 -#1)
\savepos(#1 #1)(*ex *ey) \esegment \move(*ex *ey)}

\def\dblnw#1{\bsegment
\rmove(.05 .05) \rlvec(-#1 #1) \rmove(-.1 -.1) \rlvec(#1 -#1)
\savepos(-#1 #1)(*ex *ey) \esegment \move(*ex *ey)}

\newbox\Adr
\setbox\Adr\vbox{
\vskip.5cm
\centerline{Institut für Mathematik der Universität Wien,}
\centerline{Strudlhofgasse 4, A-1090 Wien, Austria.}
\centerline{E-mail: \footnotesize{\tt Theresia.Eisenkoelbl@univie.ac.at}}
}

\title{2--enumerations of halved alternating sign matrices}
\author{Theresia Eisenkölbl\\
\box\Adr
}
\subjclass{Primary 05A15 ;
 Secondary 05C70}
\keywords{}
\begin{abstract}
We compute $2$--enumerations of certain halved alternating sign 
matrices. In one case the enumeration equals the number of perfect
matchings of a halved Aztec diamond. In the other case the enumeration 
equals the number of perfect matchings of a halved fortress graph. Our 
results prove three conjectures by Jim Propp.
\end{abstract}
\maketitle
An alternating sign matrix is a square matrix with entries $0,1,-1$
where the entries 1 and $-1$ alternate in each row and column and the
sum of entries in each row and column is equal to 1. An example of an
alternating sign matrix of order 6 is
$$\begin{pmatrix}0&0&\phantom{-}1&0& \phantom{-}0&0\\
                 1&0&-1&0& \phantom{-}1&0\\
                 0&0&\phantom{-} 1&0&-1&1\\
                 0&1&-1&1& \phantom{-}0&0\\
                 0&0&\phantom{-} 0&0& \phantom{-}1&0\\
                 0&0&\phantom{-} 1&0& \phantom{-}0&0
\end{pmatrix}.$$
Given a $k\times k$ alternating sign matrix with entries $a_{ij}$,
$1\le i,j \le k$, we form %the corner--sum
                                                     %matrix
                                                     %$(c_{ij})_{i,j=0}^n$  
the corresponding height matrix $h$ with $h_{ij}=i+j-2 \sum _{l=1}
^{i}\sum _{r=1} ^{j}a_{lr}$, $0\le i,j\le k$. 
  
The height matrix for the above alternating sign matrix looks as follows:

$$h=\begin{pmatrix}
0&1&2&3&4&5&6\\
1&2&3&2&3&4&5\\
2&1&2&3&4&3&4\\
3&2&3&2&3&4&3\\
4&3&2&3&2&3&2\\
5&4&3&4&3&2&1\\
6&5&4&3&2&1&0
\end{pmatrix}.
$$             

A height matrix has first row and column $(0,1,\dots,k)$, last row and 
column $(k,\dots,1,0)$ and adjacent entries differing by one.

Now we look at halved alternating sign matrices of order $2n$, i.e.,
$n\times 2n$--rectangles with entries $0,1, -1$ where the non--zero
entries alternate in each row and column, the row sums equal 1 and the 
first non--zero entry in each column is 1 if there is any. There is a
corresponding $(n+1)\times(2n+1)$--rectangle of heights. In this
article we only consider halved alternating sign matrices
corresponding to height matrices of the form

\begin{equation} \label{bed} 
\begin{pmatrix}
0       &1      &2      &3      &4      &\dots&2n-2     &2n-1   &2n\\
1       &?      &?      &?      &?      &\dots&?        &?      &2n-1\\
2       &?      &?      &?      &?      &\dots&?        &?      &2n-2\\
\vdots  &\vdots &\vdots &\vdots &\vdots &\dots&\vdots   &\vdots &\vdots\\
n-2     &?      &?      &?      &?      &\dots&?        &?      &n+2\\
n-1     &?      &?      &?      &?      &\dots&?        &?      &n+1\\
n       &c_1    &n      &c_2    &n      &\dots&n        &c_n    &n
\end{pmatrix}.
\end{equation}

In \cite{Propp}, Propp states conjectures regarding weighted
enumerations of halved alternating
sign matrices with height matrices of the form \eqref{bed}. Some of
these conjectures are proved in the following theorems. 
\begin{theorem} \label{dom}
The weighted enumeration with weight $2^{N_{-}(A)}$
of halved alternating sign matrices $A$ of order $2n$ with height
matrix of the form \eqref{bed}  
is $2^{n^2}$, where $N_{-}(A)$ is the number of $(-1)$'s in the halved 
alternating sign matrix $A$.
\end{theorem}

\begin{theorem} \label{fort}
The weighted enumeration with weight $2^{N_{-}(A,\,\,even)+N_{+}(A,\,\,odd)}$
of halved alternating sign matrices $A$ of order $2n$ with height
matrix of the form \eqref{bed}  
is $3^n5^{\binom n2}$, where $N_{-}(A,\,\,even)$ is the number of $(-1)$'s
in the halved  
alternating sign matrix $A$ in even position (i.e., the sum of the row index
and the column index is even) and $N_{+}(A,\,\,odd)$ is the number of 1's
in odd position.
\end{theorem}

\begin{remark}
If we use the weight $2^{N_{-}(A,\,\,odd)+N_{+}(A,\,\,even)}$, we get the same
result because reflecting a halved alternating sign matrix
corresponding to a height matrix of the form \eqref{bed} with
respect to a vertical symmetry axis gives a matrix of the same form
and interchanges even and odd 
positions of entries.
\end{remark}

\begin{remark}
The weighted enumeration of all alternating sign matrices $A$ of 
order $n$ with
weight $2^{N_{-}(A)}$ gives $2^{\binom{n}2}$, the number of perfect
matchings of an Aztec diamond of order $n-1$, see \cite{EKLP}.

The weighted enumeration of all alternating sign matrices of order
$2n$ with weight $2^{N_{-}(A,\,\,even)+N_{+}(A,\,\,odd)}$ gives
$5^{n^2}$, the number of perfect matchings in a $2n\times 2n$ fortress 
graph, see \cite[Ch.3]{yang}.
\end{remark}

\begin{theorem} \label{fix}
The weighted enumeration with weight $2^{N_{-}(A,\,\,even)+N_{+}(A,\,\,odd)}$
of halved alternating sign matrices $A$ of order $2n$ with height
matrix of the form \eqref{bed} with the additional constraint that
$c_i=n+1$ for all $i$ equals $5^{\binom n2}$ for even $n$ and
$2^n5^{\binom n2}$ for odd $n$.
If $c_i=n-1$ for all $i$, it equals $2^n5^{\binom n2}$ for even $n$
and $5^{\binom n2}$ for odd $n$.
\end{theorem} 

Below we give the proof of Theorem~\ref{dom}.
The proof of Theorem~2 starts on page~\pageref{profo}.
The proof of Theorem~3 is sketched on page~\pageref{fixproof}.

\begin{proof}[Proof of Theorem~\ref{dom}]
Since adjacent entries in the height matrix differ by one, each $c_i$ is
either $n-1$ or $n+1$.

There is a  well--known 1 to $2^{N_{-}(A)}$ correspondence between
(halved) alternating sign matrices and perfect matchings of (halved)
Aztec diamonds 
(cf. \cite{Ciucu}). An $m\times k$ Aztec rectangle is a graph composed 
of $m\times k$ squares (see Figure~\ref{azrecfi}). A halved Aztec
diamond is an Aztec rectangle with the shape of half a square, with
some vertices in the two bottom rows  missing (see Figure~\ref{newregion}).
A perfect matching (1--factor) of a graph is a set of edges such that
every vertex of the graph lies on exactly one of these edges. In the
remainder of this paper we will use the term matching instead of
perfect matching.

\begin{figure}
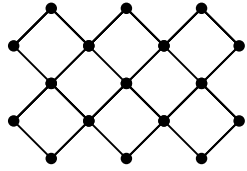

\centertexdraw{
\drawdim cm \setunitscale.5
\quadk \quadk \quadk
\move(0 -2)
\quadk \quadk \quadk
}
\caption{\label{azrecfi} A $2\times 3$ Aztec rectangle with no
vertices missing.}
\end{figure}

We write the entries of the halved alternating sign matrix in the
squares of the 
halved Aztec diamond as shown in Figure~\ref{domfi}. The corresponding 
$2^{N_{-}(A)}$ matchings can be found by demanding that a square
surrounding a $-1$, 0 or 1 contains exactly 2, 1 or 0 edges of the
matching, respectively. The edges in the squares containing a 0 can be 
found by joining the two vertices lying in the direction of the next
1's in the same row and column (if there is no 1 in the column we take the 
bottom vertex). There are two choices for the squares containing $-1$
as shown in Figure~\ref{2ways}. This accounts for the weight $2^{N_{-}(A)}$.

\begin{figure}
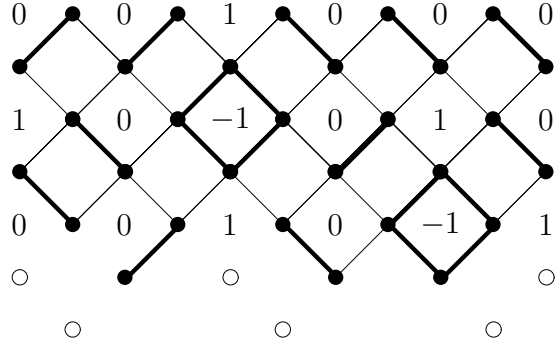

\centertexdraw{
\drawdim cm \setunitscale.7

\linewd.01
\move(0 0) \quadk \quadk \quadk \quadk \quadk
\move(0 -2) \quadk \quadk \quadk \quadk \quadk
%\move(1 -3) \quadk
%\move(5 -3) \quadk \rmove(1 -1) \knoten
%\move(0 -4) \quadi
%\move(4 -4) \quadi
\move(6 -4) \rlvec(1 1)
%\move(8 -4) \quadi
%\move(4 -4) \rlvec(-1 1)
\move(2 -4) \knoten
\move(6 -4) \knoten
\move(8 -4) \knoten

\move(0 -4)\lcir r:.15
\move(4 -4)\lcir r:.15
\move(10 -4)\lcir r:.15
\move(1 -5)\lcir r:.15
\move(5 -5)\lcir r:.15
\move(9 -5)\lcir r:.15

\linewd.08
\move(0 0)
\rlvec(1 1) \rmove(1 -1) \rlvec(1 1)  \rmove(1 -1)  \rlvec(1 -1)
\rlvec(-1 -1)  \rlvec(-1 1)  \rlvec(1 1) \rmove(1 1)  \rlvec(1 -1)
\rmove(1 1)  \rlvec(1 -1) \rmove(1 1)  \rlvec(1 -1) \rmove(-1 -1)
\rlvec(1 -1)  \rmove(-1 -1)  \rlvec(-1 -1)  \rlvec(-1 1) \rlvec(1 1)
\rlvec(1 -1)  \rmove(-2 2)  \rlvec(-1 -1)  \rmove(0 -2)  \rlvec(-1 1)
\rmove(-2 0)  \rlvec(-1 -1) \rmove(-1 1)  \rlvec(-1 1) \rmove(1 1)
\rlvec(1 -1)  

\linewd.05
%\move(9 -5) \dblne1
%\move(1 -5) \dblnw1
%\move(5 -5) \dblnw1
\textref h:C v:C
\htext(0 1){$0$}
\htext(2 1){$0$}
\htext(4 1){$1$}
\htext(6 1){$0$}
\htext(8 1){$0$}
\htext(10 1){$0$}
\htext(0 -1){$1$}
\htext(2 -1){$0$}
\htext(4 -1){$-1$}
\htext(6 -1){$0$}
\htext(8 -1){$1$}
\htext(10 -1){$0$}
\htext(0 -3){$0$}
\htext(2 -3){$0$}
\htext(4 -3){$1$}
\htext(6 -3){$0$}
\htext(8 -3){$-1$}
\htext(10 -3){$1$}

}
\caption{\label{domfi} The set of perfect matchings corresponding to a
halved alternating sign matrix. }
\end{figure}

\begin{figure}
\centertexdraw{
\drawdim cm \setunitscale.7
\linewd.02 
\move(4 0)\quadk
\move(8 0)\quadk
\linewd.1
\move(0 0)
\quadk \move(4 0) \rlvec(1 1) \move(5 -1) \rlvec(1 1) 
\move(8 0) \rlvec(1 -1) \move(9 1) \rlvec(1 -1) 
\textref h:C v:C
\htext(3 0){=}
\htext(7 0){+}

}
\caption{\label{2ways} }
\end{figure}

Close inspection of the correspondence reveals that the condition on
the last row of the height matrix determines which of the vertices of
the halved Aztec diamond are missing.
The halved Aztec diamond is an $n\times (2n-1)$ Aztec
rectangle with missing last row of vertices and missing vertices in the
next row in positions $a_1, \dots, a_n$, say. It is easy to see that we
have either $a_i=2i-1$ or $a_i=2i$ corresponding to $c_i=n-1$ or
$n+1$, respectively. 
Therefore, we have to sum
over $2^n$ different boundary conditions. Fortunately, we can add
pairs of vertices in the last two rows as shown in Figure~\ref{domfi}
and just count all matchings of the emerging new region (see 
Figure~\ref{newregion}). The
vertices in the bottom row can be matched either to the northeast or to
the northwest. This corresponds to the possible choices for the $a_i$.

\begin{figure}
\centertexdraw{
\drawdim truecm \setunitscale.7
\linewd.05
\quadk\quadk\quadk\quadk\quadk
\move(0 -2) \quadk \quadk\quadk\quadk\quadk
\move(0 -4) \quadk \rlvec(1 1)\rlvec(1 -1)\quadk\rlvec(1 1) \rlvec(1
-1)\quadk

}
\caption{\label{newregion} }
\end{figure}

Now we can apply the following lemma (cf. \cite[p.18]{EKLP}).
\begin{lemma} \label{azrec}
The number of perfect matchings of an $m\times k$ Aztec rectangle, where
all the vertices in the bottom row have been removed except for the
$x_1$st, the $x_2$nd, \dots, and the $x_m$th vertex
equals
$$\frac {2^{\binom {m+1}2}} {\prod _{i=1} ^{m}(i-1)!}\prod _{1\le
i<j\le m} (x_j-x_i).$$ 
\end{lemma}

To apply the lemma to our case, we have to set $k=2n-1$, $m=n$ and
$x_i=2i-1$. We obtain that our 2--enumeration of halved alternating
sign matrices equals
$$\frac {2^{\binom {n+1}2}} {\prod _{i=1} ^{n}(i-1)!}\prod _{1\le
i<j\le n} (2j-2i)=2^{\binom{n+1}2}2^{\binom n2}=2^{n^2},$$ 
as desired.
\end{proof}

\begin{proof}[Proof of Theorem~\ref{fort}] \label{profo}
Now we have the weight
$2^{N_{-}(A,\,\,even)+N_{+}(A,\,\,odd)}$. 
We will illustrate all steps of the proof by the example of halved
alternating sign matrices of order 6 (cf. Figure~\ref{fortfi}). From
there it will be clear what happens in the general case.

The first step is another well--known
bijection between alternating sign matrices and a family of graphs 
called fortresses. These are squares arranged in a rectangular shape
separated by single edges. On the right, the left and the upper side
of the rectangle edges are appended to every other square
(see Figure~\ref{fortfi} for a
$3\times6$ fortress graph with some extra edges appended to the
squares in the bottom row). 
We have the following replacement rules:
\begin{itemize}
\item 1's in {\em even} places and $-1$'s in {\em odd} places
translate to
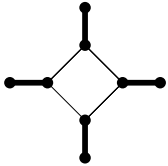
\nopagebreak

 \begin{texdraw} \drawdim cm \setunitscale.5
\linewd.01 \urba \linewd.15 \rlvec(-1 0) \rmove(-1 1) \rlvec(0 1)
\move(0 0) \rlvec(1 0) \rmove(1 -1) \rlvec(0 -1)
\end{texdraw}

\item 
$-1$'s in {\em even} places and 1's in {\em odd} places
translate to
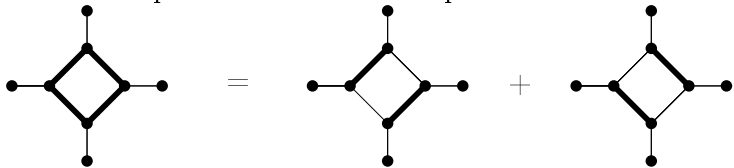
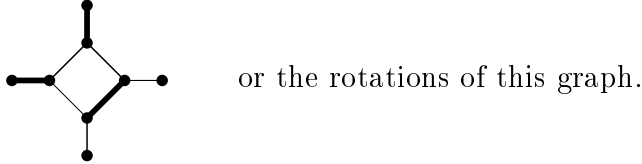
\nopagebreak

\begin{texdraw} \drawdim cm \setunitscale.5
\linewd.01 \urba \linewd.15 \move(1 0) \quadi
\textref h:C v:C
\htext(6 0){\boldmath =}
\move(8 0)\linewd.01 \urba \linewd.15 \rmove(-3 0) \rlvec(1 1) \rmove(1
-1) \rlvec(-1 -1)
\htext(13.5 0){\boldmath +}
\move(15 0)
\linewd.01 \urba \linewd.15 \rmove(-3 0)  \rlvec(1 -1)
\rmove(1 1) \rlvec(-1 1)
\end{texdraw}

\item 
0's translate to
\nopagebreak

\begin{texdraw} \drawdim cm \setunitscale.5
\linewd.01 \urba \linewd.15 \move(0 0) \rlvec(1 0) \rmove(1 -1) \rlvec(1 1)
\rmove(-1 1) \rlvec(0 1)
\textref h:L v:C
\htext(6 0){or the rotations of this graph.}
\end{texdraw}
\end{itemize}

The edges of the squares corresponding to $\pm 1$ determine uniquely
which of the four possibilities should be chosen for each 0.
An example of the correspondence is shown in Figure~\ref{fortfi}. The
reader should note that there are two choices of edges for 
$-1$'s in even places and for 1's in odd places. This accounts for the
weight. 

\begin{figure}
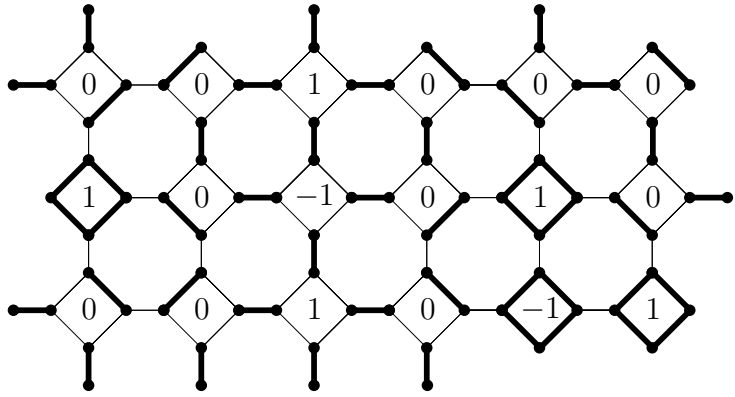

\centertexdraw{
\drawdim cm \setunitscale.5
\linewd.02
\urba\quadk\urba\quadk\urba\quadk
\move(1 -3)\quadk\urba\quadk\urba\quadk\urba
\move(0 -6)\urba\quadk\urba\quadk\rlvec(1 0) \knoten 
\quadk \rmove(-1 1) \rlvec(0 1) \rmove(1 -2) \rlvec(1 0)\quadk
\move(5 -7) \rlvec(0 -1) \knoten
\move(11 -7) \rlvec(0 -1)\knoten
%\move(14 -7) \rlvec(0 -1)\knoten
%\move(17 -7) \rlvec(0 -1)\knoten
\linewd.15
\move(0 0) \rlvec(1 0)
\rmove(1 1) \rlvec(0 1)
\move(0 -6) \rlvec(1 0)
\move(2 -1) \rlvec(1 1)
\move(1 -3) \quadi
\move(2 -5) \rlvec(1 -1) \rmove(-1 -1) \rlvec(0 -1)
\move(4 0) \rlvec(1 1) \rmove(1 -1) \rlvec(1 0)
\rmove(-2 -1) \rlvec(0 -1) \rmove(-1 -1) \rlvec(1 -1) \rmove(1 1)
\rlvec(1 0) \rmove(-2 -2) \rlvec(-1 -1) \rmove(1 -1) \rlvec(0 -1)
\rmove(1 2) \rlvec(1 0) \rmove(1 -1) \rlvec(0 -1) \rmove(1 2) 
\rlvec(1 0) \rmove(0 3) \rlvec(-1 0)
\move(8 2) \rlvec(0 -1) \rmove(0 -2) \rlvec(0 -1) \rmove(0 -2)
\rlvec(0 -1)
\move(9 0) \rlvec(1 0) \rmove(1 1) \rlvec(1 -1) \rmove(-1 -1) 
\rlvec(0 -1) 
\rmove(1 -1) \rlvec(-1 -1) \rmove(0 -1)  \rlvec(1 -1) \rmove(-1 -1)
\rlvec(0 -1)
\move(14 2) \rlvec(0 -1) \rmove(-1 -1) \rlvec(1 -1) \rmove(0 -1)
\rmove(-1 -1) \quadi \rmove(-2 -3) \quadi \rmove(1 0) 
\quadi \rmove(-2 3) \rlvec(1 -1)
\rmove(1 1) \rlvec(1 0) \rmove(-2 1) \rlvec(0 1) \rmove(1 1) \rlvec(-1
1) \rmove(-1 -1) \rlvec(-1 0)
\linewd.05
\textref h:C v:C

\htext(2 0){$0$}
\htext(5 0){$0$}
\htext(8 0){$1$}
\htext(11 0){$0$}
\htext(14 0){$0$}
\htext(17 0){$0$}
\htext(2 -3){$1$}
\htext(5 -3){$0$}
\htext(8 -3){$-1$}
\htext(11 -3){$0$}
\htext(14 -3){$1$}
\htext(17 -3){$0$}
\htext(2 -6){$0$}
\htext(5 -6){$0$}
\htext(8 -6){$1$}
\htext(11 -6){$0$}
\htext(14 -6){$-1$}
\htext(17 -6){$1$}
}
\caption{\label{fortfi} The corresponding matchings of the fortress graph. }
\end{figure}

It is not difficult to see that the restriction on the last row of the 
height matrix corresponds to a condition on the extra pending edges at 
the bottom row of the resulting graph. Either both edges in the
positions $2i$ and $2i-1$ are contained in the graph or neither.

\begin{figure}
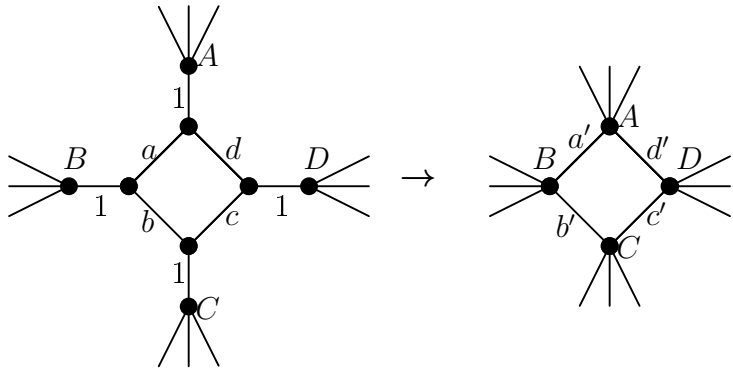

\centertexdraw{
\drawdim cm
\setunitscale.8
 \quadk
\ringerl(0 0) \ringerl(2 0) \ringerl(1 1) \ringerl(1 -1)
\ringerl(-1 0) \ringerl(3 0) \ringerl(1 2) \ringerl(1 -2)
\move(2 0)  \rlvec(1 0)
\rlvec(1 .5) \move(3 0) \rlvec(1 0) \move(3 0) \rlvec(1 -.5)
  \move(1 1) \rlvec(0 1) 
   \rlvec(-.5 1) \move(1 2) \rlvec(0 1) \move(1 2) \rlvec(.5 1)
  \move(0 0) \rlvec(-1 0)
   \rlvec(-1 .5) \move(-1 0) \rlvec(-1 0) \move(-1 0) \rlvec(-1 -.5)
  \move(1 -1) \rlvec(0 -1) 
 \rlvec(.5 -1) \move(1 -2) \rlvec(0 -1) \move(1 -2) \rlvec(-.5 -1) 
 \htext(.2 .45){$a$}
 \htext(1.6 .45){$d$}
 \htext(1.6 -.65){$c$}
 \htext(.2 -.75){$b$}
\htext(1.1 2){$A$}
\htext(1.1 -2.2){$C$}
\htext(-1.1 .3){$B$}
\htext(2.9 .3){$D$}
\htext(2.4 -.5){1}
\htext(-.6 -.5){1}
\htext(.7 1.3){1}
\htext(.7 -1.6){1}

\move(7 0)
\bsegment 
  \quadk
\ringerl(0 0) \ringerl(2 0) \ringerl(1 1) \ringerl(1 -1)
\move(2 0)
   \rlvec(1 .5) \move(2 0) \rlvec(1 0) \move(2 0) \rlvec(1 -.5)
  \move(1 1) 
   \rlvec(-.5 1) \move(1 1) \rlvec(0 1) \move(1 1) \rlvec(.5 1)
  \move(0 0) 
   \rlvec(-1 .5) \move(0 0) \rlvec(-1 0) \move(0 0) \rlvec(-1 -.5)
  \move(1 -1)  
\rlvec(.5 -1) \move(1 -1) \rlvec(0 -1) \move(1 -1) \rlvec(-.5 -1) 
 \htext(.3 .65){$a'$}
 \htext(1.6 .45){$d'$}
 \htext(1.6 -.65){$c'$}
 \htext(.1 -.85){$b'$}
\htext(1.1 1){$A$}
\htext(1.1 -1.2){$C$}
\htext(-0.3 .3){$B$}
\htext(2.1 .3){$D$}
\esegment
\htext(4.5 0){\boldmath $\rightarrow$}
}
\caption{\label{urbfi} Urban Renewal, $a'=\frac {c} {ac+bd}$, $b'=\frac
{d} {ac+bd}$, $c'=\frac {a} {ac+bd}$, $d'=\frac{b}{ac+bd}$.}
\end{figure}

In the following we will repeatedly use a well--known local
modification of a graph called
urban renewal, which changes the enumeration of perfect
matchings only by a global factor (see \cite{shuffle}). The
modification is shown in Figure~\ref{urbfi}.
Before we can
explain this modification, we have to make a few definitions.
Let $G$ be a graph with weights assigned to its edges. Then the weight 
of a matching is the product of the weights of the edges it contains.
The weighted enumeration of matchings $M(G)$ is now defined as the sum 
of the weights of all possible matchings of the graph $G$.

Urban renewal can now be described as follows. We start with a graph
$G$ which looks locally like the left--hand--side of
Figure~\ref{urbfi}. Then we contract the four edges of weight 1
and change the weights $a,b,c,d$ to $a',b',c',d'$. We obtain a graph
$G'$ which looks locally like the right--hand--side of
Figure~\ref{urbfi} and like $G$ everywhere else.
The new edge weights $a',b',c',d'$ of the resulting $G'$ are defined by
$a'=\frac {c} {ac+bd}, b'=\frac {d} {ac+bd},
c'=\frac {a} {ac+bd}, d'=\frac{b}{ac+bd}$, whereas all other weights
stay the same. 

The weighted enumerations of matchings $M(G)$ and $M(G')$ of the two
graphs are related in the following way:

\begin{lemma} \label{urb}
Let $G$ be a graph which looks locally like the left--hand--side of
Figure~\ref{urbfi} and let $G'$ be the graph which looks locally like
the right--hand--side and like $G$ elsewhere.
Then the weighted enumeration of matchings of the new graph
$G'$ equals the weighted enumeration of matchings of the graph $G$ 
multiplied by $ac+bd$, i.e.,
$$M(G)=M(G')(ac+bd).$$
\end{lemma}

\begin{figure}
\centertexdraw{\drawdim cm \setunitscale.5
\ringquad\quadk\ringquad\quadk\ringquad\quadk
\move(1 -3)\quadk\ringquad\quadk\ringquad\quadk\ringquad
\move(0 -6)\ringquad\quadk\ringquad\quadk\ringquad\quadk
\move(5 -7) \rlvec(0 -1) \knoten
\move(11 -7) \rlvec(0 -1)\knoten
\move(14 -8) \rlvec(0 -1)\knoten
\move(0 -16)
\bsegment \setunitscale.7
\squadk\quadk\squadk\quadk\squadk\quadk
\move(0 -2) \quadk\squadk\quadk\squadk\quadk\squadk
\move(0 -4)\squadk\quadk\squadk\quadk\squadk\quadk
\move(3 -5) \rlvec(0 -1) \knoten
\move(7 -5) \rlvec(0 -1) \knoten
\move(9 -5) \rlvec(0 -1) \knoten
\esegment
\htext(8 -13){{\huge $\downarrow$} {\large $2^{n^2}$ }}
}
\caption{\label{fort2az}The dotted lines have weight $\frac {1} {2}$.}
\end{figure}

We want to apply urban renewal to all squares in even position in the
graph in Figure~\ref{fortfi}.
First, we append two vertical edges to squares in the last row in even
position which have no downward--pointing edge appended. In the example 
in Figure~\ref{fortfi}, this happens to the fifth square in the bottom 
row, resulting in the upper graph in Figure~\ref{fort2az} (at this
point, the dotted circles should be ignored).
This does not 
change the enumeration of perfect matchings since there is only one
possibility for the new vertices to be paired.
We obtain an
$n\times 2n$ fortress with some edges appended.
% (see the upper graph in Figure~\ref{fort2az}).
Now, we can apply urban renewal to every square in even position (the
circled squares in Figure~\ref{fort2az}).  
There are $n^2$ of these squares with $ac+bd=2$,
which yields a factor of $2^{n^2}$.
The resulting graph is an $n\times 2n$ Aztec rectangle where every
other square has edges of weight $\frac {1} {2}$ (see the dotted lines
in the bottom graph
in Figure~\ref{fort2az}) and some downward--pointing edges appended to 
the last row of squares. It is easy to see that the original
restriction on the last row of the height matrix now
translates to the restriction that for each $i$
there is exactly one edge in the
position $2i-1$ or $2i$. Similar to the proof of Theorem~\ref{dom}, 
we can interpret the sum of the corresponding $2^n$ terms
as the number of matchings of the weighted halved Aztec diamond $G_n$ in
Figure~\ref{sumaz} because every vertex in its last row can be either
matched to the left or to the right.

Therefore, we have to determine $2^{n^2}\times M(G_n)$.

\begin{figure}
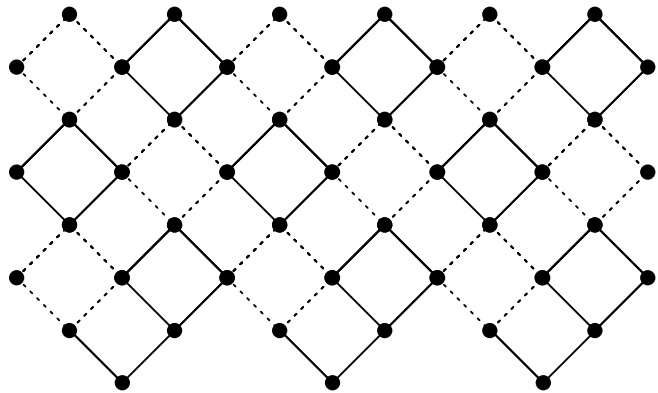

\centertexdraw{
\drawdim cm
\setunitscale.7
\squadk\quadk\squadk\quadk\squadk\quadk
\move(0 -2) \quadk\squadk\quadk\squadk\quadk\squadk
\move(0 -4)\squadk\quadk\squadk\quadk\squadk\quadk
\move(1 -5) \rlvec(1 -1) \knoten \rlvec(1 1)
\move(5 -5) \rlvec(1 -1) \knoten \rlvec(1 1)
\move(9 -5) \rlvec(1 -1) \knoten \rlvec(1 1)
}
\caption{\label{sumaz} $G_n$ %, the sum of the $2^n$ possible
                             %configurations 
for $n=3$.}
\end{figure}

We will use urban renewal repeatedly to reduce the graph $G_n$ to the
graph $G_{n-1}$. The first step consists of replacing every vertex of
$G_n$ by a line of three vertices so that the weighted enumeration of
matchings remains unchanged. The resulting graph in our example is shown in
Figure~\ref{split3}. Now we are in the position to apply urban renewal
to all squares in the graph. The factor $ac+bd$ equals $\frac {1}
{2}\cdot \frac 12 +\frac {1} {2}\cdot\frac {1} {2}=\frac {1} {2}$ for
half of the 
squares and $1\cdot1+1\cdot1=2$ for the other half. Thus, the factors
resulting from urban renewal cancel each other. Square edges of 
weight 1 become edges of weight $\frac {1} {2}$ and vice versa and in
our example we
obtain the graph shown in Figure~\ref{trafo}.
\begin{figure}
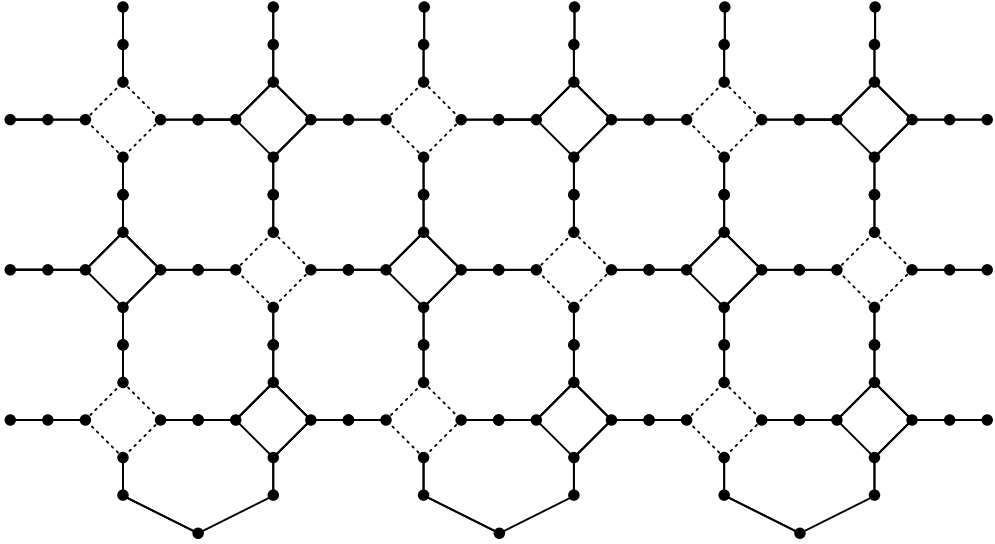

\centertexdraw{
\drawdim truecm \setunitscale.5
\knoten \rlvec(1 0) \urbb \urba \urbb \urba \urbb \urba \rlvec(1 0)
\knoten
\move(0 -4)
\knoten \rlvec(1 0)  \urba \urbb \urba \urbb \urba \urbb \rlvec(1 0) \knoten
\move(0 -8) 
\knoten \rlvec(1 0) \urbb \urba \urbb \urba \urbb \urba \rlvec(1 0)
\knoten
\move(3 -10) \rlvec(2 -1) \knoten \rlvec(2 1) 
\move(11 -10) \rlvec(2 -1) \knoten \rlvec(2 1) 
\move(19 -10) \rlvec(2 -1) \knoten \rlvec(2 1) 
\move(3 2) \rlvec(0 1) \knoten
\move(7 2) \rlvec(0 1) \knoten
\move(11 2) \rlvec(0 1) \knoten
\move(15 2) \rlvec(0 1) \knoten
\move(19 2) \rlvec(0 1) \knoten
\move(23 2) \rlvec(0 1) \knoten

}
\caption{\label{split3} Each vertex is split into three.}
\end{figure}

\begin{figure}
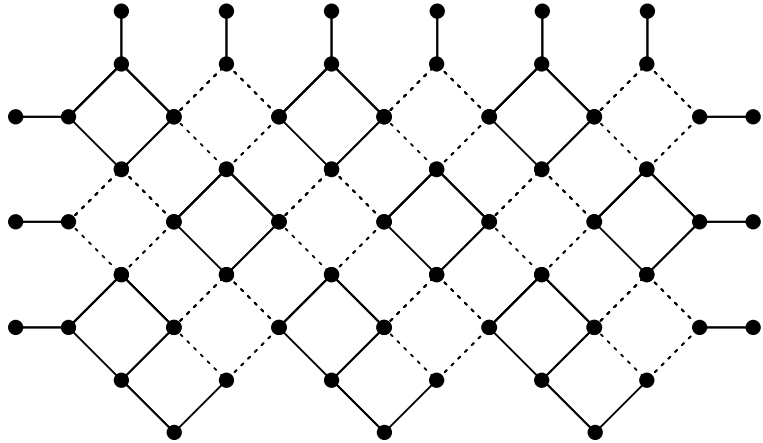

\centertexdraw{
\drawdim cm \setunitscale.7
\quadk\squadk\quadk\squadk\quadk\squadk
\move(0 -2)
\squadk\quadk\squadk\quadk\squadk\quadk
\move(0 -4)
\quadk\squadk\quadk\squadk\quadk\squadk
\move(0 0) \rlvec(-1 0) \knoten
\move(0 -2) \rlvec(-1 0) \knoten
\move(0 -4) \rlvec(-1 0) \knoten
\move(1 1) \rlvec(0 1) \knoten
\move(3 1) \rlvec(0 1) \knoten
\move(5 1) \rlvec(0 1) \knoten
\move(7 1) \rlvec(0 1) \knoten
\move(9 1) \rlvec(0 1) \knoten
\move(11 1) \rlvec(0 1) \knoten
\move(12 0) \rlvec(1 0) \knoten
\move(12 -2) \rlvec(1 0) \knoten
\move(12 -4) \rlvec(1 0) \knoten
\move(1 -5) \rlvec(1 -1) \knoten \rlvec(1 1)
\move(5 -5) \rlvec(1 -1) \knoten \rlvec(1 1)
\move(9 -5) \rlvec(1 -1) \knoten \rlvec(1 1)

}
\caption{\label{trafo} The graph obtained by applying urban renewal to 
all squares in Figure~\ref{split3}.}
\end{figure}

The pending edges along the border of the graph have to be in every perfect
matching and can be removed together with the two endpoints without
changing the enumeration of perfect matchings. For the same reason, we 
can fill the ``dents'' in the bottom row by adding some edges which
have to be in every perfect matching. The resulting graph is shown in
Figure~\ref{+-fi}. 

\begin{figure}
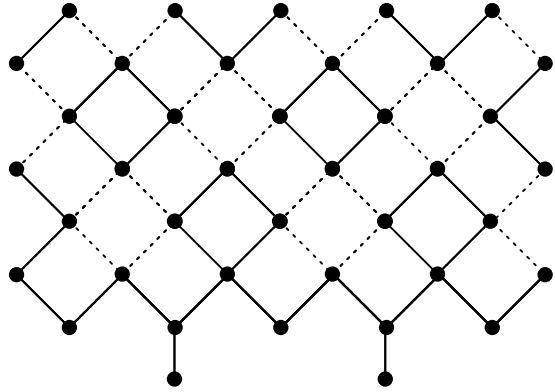

\centertexdraw{
\drawdim cm \setunitscale.7
\rlvec(1 1) \knoten \lpatt(.05 .13) \rlvec(1 -1) 
\rlvec(1 1) \knoten \lpatt(1 0) \rlvec(1 -1) 
\rlvec(1 1) \knoten \lpatt(.05 .13) \rlvec(1 -1)
\rlvec(1 1) \knoten \lpatt(1 0) \rlvec(1 -1) 
\rlvec(1 1) \knoten \lpatt(.05 .13) \rlvec(1 -1) \lpatt(1 0) \knoten 
\rlvec(-1 -1) \rlvec(1 -1) \knoten  \lpatt(.05 .13) \rlvec(-1 -1) \rlvec(1 -1)
\lpatt(1 0) \knoten 
\rlvec(-1 -1) \knoten  \rlvec(-1 1)
\rlvec(-1 -1) \knoten \rlvec(-1 1)
\rlvec(-1 -1) \knoten \rlvec(-1 1)
\rlvec(-1 -1) \knoten \rlvec(-1 1)
\rlvec(-1 -1) \knoten  \rlvec(-1 1) \knoten
\rlvec(1 1) \rlvec(-1 1) \knoten
\lpatt(.05 .13) 
\rlvec(1 1) \rlvec(-1 1) \knoten \lpatt(1 0)
\move(3 -5) \rlvec(0 -1) \knoten
\move(7 -5) \rlvec(0 -1) \knoten

\move(1 -1) \quadk \squadk \quadk \squadk
\move(1 -3) \squadk \quadk \squadk \quadk
}
\caption{\label{+-fi} Remove and add some forced edges of weight 1.}
\end{figure}

\begin{figure}
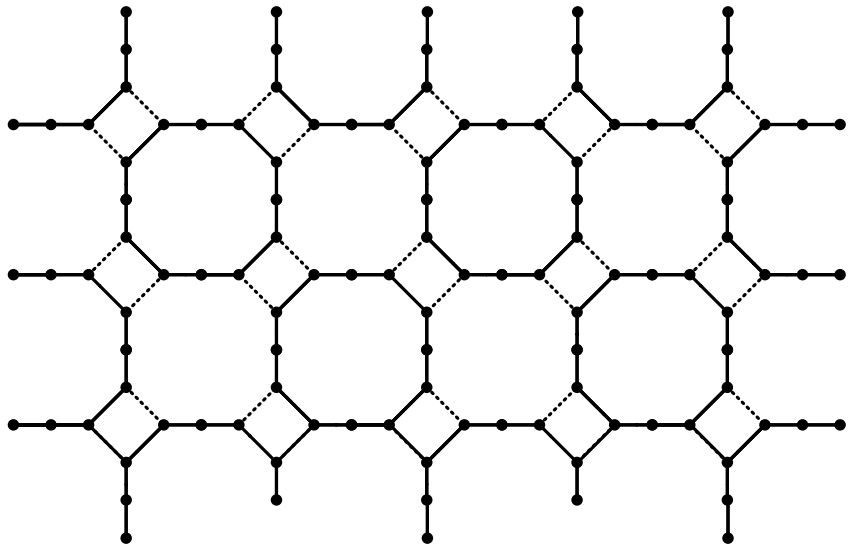

\centertexdraw{
\drawdim cm \setunitscale.5 \linewd.08
\knoten \rlvec(1 0) \urbd \urbc \urbd \urbc \urbd  \rlvec(1 0)
\knoten
\move(0 -4)
\knoten \rlvec(1 0)  \urbc \urbd \urbc \urbd \urbc \rlvec(1 0) \knoten
\move(0 -8) 
\knoten \rlvec(1 0) \urbd \urbc \urbd \urbc \urbd  \rlvec(1 0)
\knoten
\move(3 -10) \rlvec(0 -1) \knoten  
\move(11 -10) \rlvec(0 -1) \knoten  
\move(19 -10) \rlvec(0 -1) \knoten 
\move(3 2) \rlvec(0 1) \knoten
\move(7 2) \rlvec(0 1) \knoten
\move(11 2) \rlvec(0 1) \knoten
\move(15 2) \rlvec(0 1) \knoten
\move(19 2) \rlvec(0 1) \knoten

\move(2 -8) \rlvec(1 -1)
\move(7 -9) \rlvec(1 1)
\move(10 -8) \rlvec(1 -1)
\move(15 -9) \rlvec(1 1)
\move(18 -8) \rlvec(1 -1)
}
\caption{\label{split3bfi} The graph obtained by replacing all vertices in
Figure~\ref{+-fi} by three vertices.}
\end{figure}

Now we (almost) repeat the last two steps. We replace 
each vertex by three
vertices to obtain the graph in Figure~\ref{split3bfi}. Note that the
squares in the
bottom row contain only one edge of weight $\frac {1} {2}$ each.

The next step is to apply urban renewal to all squares. The product of
the factors $ac+bd$ is
easily seen to be $\(\frac {5} {4}\)^{(n-1)(2n-1)}\(\frac {3}
{2}\)^{2n-1}$.
The new edge weights are $\frac {2} {3}$, $\frac {1} {3}$, $\frac {2}
{5}$ and $\frac {4} {5}$ (only the forced appended edges have still
weight 1). The resulting graph is shown in Figure~\ref{lowfi}.

Now we mark every other vertex in the bottommost row with a dotted
circle starting with the second vertex (see Figure~\ref{lowfi}). 
Similarly, we mark all vertices immediately above and to the left of
the dotted circles with an unbroken circle.
We divide the weight of the edges incident to one of the $n-1$
points marked by an unbroken circle by two and
multiply the 
weight of the edges incident to the $n-1$ points marked by a dotted
circle by two. This does not change the weighted
enumeration of matchings.

\begin{figure}
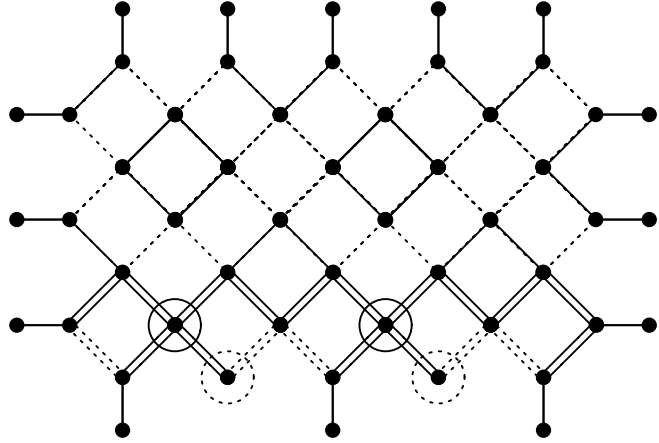

\centertexdraw{
\drawdim cm \setunitscale.7
\squadk\squadk\squadk\squadk\squadk
\move(0 -2) \squadk \squadk\squadk \squadk \squadk
\move(1 -1) \quadk \squadk \quadk \squadk
%\move(1 -3) \quadk \squadk \quadk \squadk
\move(1 1) \rlvec(0 1) \knoten
\move(3 1) \rlvec(0 1) \knoten
\move(5 1) \rlvec(0 1) \knoten
\move(7 1) \rlvec(0 1) \knoten
\move(9 1) \rlvec(0 1) \knoten
\move(10 0) \rlvec(1 0) \knoten
\move(10 -2) \rlvec(1 0) \knoten
\move(10 -4) \rlvec(1 0) \knoten
\move(0 0) \rlvec(-1 0) \knoten
\move(0 -2) \rlvec(-1 0) \knoten
\move(0 -4) \rlvec(-1 0) \knoten
\move(0 0) \rlvec(1 1) \rmove(2 0) \rlvec(1 -1) \rlvec(1 1) \rmove(2
0) \rlvec(1 -1) \rlvec(1 1) \rmove(1 -1) \rlvec(-1 -1) \rlvec(1 -1)
\move(0 -2) \rlvec(1 -1) \rmove(2 0) \rlvec(1 1) \rlvec(1 -1) \rmove(2 
0) \rlvec(1 1) \rlvec(1 -1)
\move(0 -4) \knoten \dblne1 \knoten 
\move(2 -4) \knoten \dblnw1
\move(2 -4) \knoten \dblne1 \knoten 
\move(4 -4) \knoten \dblnw1
\move(4 -4) \knoten \dblne1 \knoten 
\move(6 -4) \knoten \dblnw1
\move(6 -4) \knoten \dblne1 \knoten 
\move(8 -4) \knoten \dblnw1
\move(8 -4) \knoten \dblne1 \knoten 
\move(10 -4) \knoten \dblnw1

\move(1 -5) \knoten \lpatt(.05 .13) \dblnw1 \lpatt(1 0)
\move(1 -5) \knoten \dblne1 \knoten 
\move(3 -5) \knoten \dblnw1
\move(3 -5) \knoten \lpatt(.05 .13)\dblne1 \knoten 
\move(5 -5) \knoten \dblnw1 \lpatt(1 0)
\move(5 -5) \knoten \dblne1 \knoten 
\move(7 -5) \knoten \dblnw1
\move(7 -5) \knoten \lpatt(.05 .13)\dblne1 \knoten 
\move(9 -5) \knoten \dblnw1 \lpatt(1 0)
\move(9 -5) \knoten \dblne1 \knoten 

\move(1 -5) \rlvec(0 -1) \knoten 
\move(5 -5) \rlvec(0 -1) \knoten 
\move(9 -5) \rlvec(0 -1) \knoten 

\move(2 -4) \lcir r:.5
\move(6 -4) \lcir r:.5
\move(3 -5) \lpatt(.05 .13) \lcir r:.5
\move(7 -5) \lpatt(.05 .13) \lcir r:.5

}
\caption{\label{lowfi} The double edges have weight $\frac 23$, the
dotted double edges have weight $\frac {1} {3}$, the dotted edges have 
weight $\frac 25$, the pending edges have weight 1 and the remaining
edges have weight $\frac {4} {5}$. %At the points with the large
%circles we halve all the edges incident to them. At the points with
%the large dotted cycles we double all edges incident to them. 
}
\end{figure}

\begin{figure}
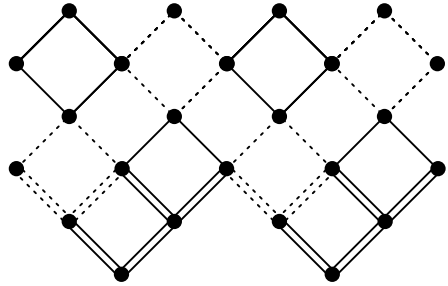

\centertexdraw{
\drawdim cm \setunitscale.7
\quadk \squadk \quadk \squadk
%\move(0 -2)
%\squadk \quadk \squadk \quadk
\move(2 -4) \dblnw1 
\move(2 -4) \knoten \dblne1
\move(6 -4) \dblnw1 
\move(6 -4) \knoten \dblne1
\move(3 -3) \dblnw1 
\move(3 -3) \knoten \dblne1
\move(7 -3) \dblnw1 
\move(7 -3) \knoten \dblne1 \knoten
\move(2 -2) \rlvec(1 1)
\move(4 -2) \rlvec(-1 1)
\move(6 -2) \rlvec(1 1)
\move(8 -2) \rlvec(-1 1) 
\lpatt(.05 .13)
\move(1 -3) \dblnw1 \knoten
\move(1 -3) \knoten \dblne1 \knoten
\move(5 -3) \dblnw1 \knoten
\move(5 -3) \knoten \dblne1 \knoten
\move(0 -2) \rlvec(1 1) 
\move(2 -2) \rlvec(-1 1)
\move(4 -2) \rlvec(1 1)
\move(6 -2) \rlvec(-1 1) 

}
\caption{\label{stripfi} 
The double edges have weight $\frac 23$, the
dotted double edges have weight $\frac {1} {3}$, the dotted edges have 
weight $\frac 25$ and the remaining edges have weight $\frac {4} {5}$.}
\end{figure}

Then we strip off all the forced edges (i.e., edges that must be
contained in {\em every} perfect matching) and obtain the
graph shown in Figure~\ref{stripfi}.
It is easy to see that every matching contains exactly $2n-2$ of the
edges with weight $\frac {1} {3}$ or $\frac {2} {3}$ (the
double edges and dotted double edges) and exactly
%$((2n-1)(n-1)+(2n-2)(n-1)-(n-1))/2$ 
$2(n-1)^2$ edges with weight $\frac 25$ or $\frac 45$.
If we now divide the
weights of all double edges and dotted double edges
by $\frac {2} {3}$ and the weights of all
the other edges by $\frac {4} {5}$, we obtain a graph with edges of
weight 1 and $\frac {1} {2}$ only. This changes the weighted
enumeration by a factor of $\(\frac {2}
{3}\)^{2n-2}\(\frac {4} {5}\)^{2(n-1)^2}$.
 
The resulting graph is shown in Figure~\ref{indfi}. It  is clearly the 
mirror image of
$G_{n-1}$ (compare Figure~\ref{sumaz}). 
\begin{figure}
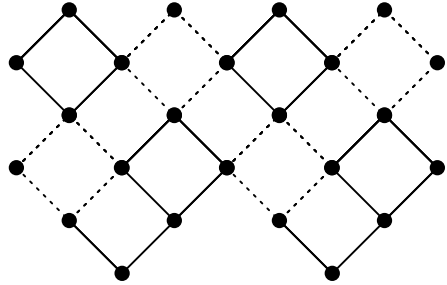

\centertexdraw{
\drawdim cm \setunitscale.7
\quadk \squadk \quadk \squadk
\move(0 -2)
\squadk \quadk \squadk \quadk
\move(1 -3) \rlvec(1 -1) \knoten \rlvec(1 1)
\move(5 -3) \rlvec(1 -1) \knoten \rlvec(1 1)
}
\caption{\label{indfi} The dotted edges have weight $\frac {1} {2}$,
the other edges have weight 1. }
\end{figure}

Therefore, we obtain for the weighted enumeration of matchings of
$G_n$:
 
$$M(G_n)=\(\frac {5} {4}\)^{(n-1)(2n-1)}\(\frac {3}
{2}\)^{2n-1}\(\frac {2}
{3}\)^{2n-2}\(\frac {4} {5}\)^{2(n-1)^2} M(G_{n-1})=\frac{3 \cdot
5^{n-1}}{2^{2n-1}} M(G_{n-1}).$$

Since $M(G_1)$ is easily seen to be $\frac {3} {2}$, we get for our
weighted enumeration of halved alternating sign matrices
$$2^{n^2}M(G_n)=2^{n^2} \frac{3^n 5^{\binom n2}}{2^{n^2}}=3^n5^{\binom
n2}.$$ 
\end{proof}

\begin{proof}[Sketch of the proof of Theorem~\ref{fix}]
\label{fixproof}
The proof is analogous to the proof of Theorem~\ref{fort}.
For example, in the case $c_i=n-1$ for all $i$, we use the bijection to
fortress graphs and apply urban renewal to all the squares. We obtain
a weighted halved Aztec diamond of order $2n$ 
(see Figure~\ref{n-1fortfi} for the case $n=3$).
Imitating the steps in
the proof of Theorem~\ref{fort}, we can reduce it to the weighted
halved Aztec diamond
of order $2n-2$. In this way, we again obtain a simple
recursion which gives the results stated in Theorem~\ref{fix}.
\begin{figure}
\centertexdraw{
\drawdim cm 
\setunitscale.7
\squadk\quadk\squadk\quadk\squadk\quadk
\move(0 -2) \quadk\squadk\quadk\squadk\quadk\squadk
\move(0 -4)\squadk\quadk\squadk\quadk\squadk\quadk
\move(3 -5) \rlvec(0 -1) \knoten
\move(7 -5) \rlvec(0 -1) \knoten
\move(11 -5) \rlvec(0 -1) \knoten
}
\caption{\label{n-1fortfi} }
\end{figure}
\end{proof}


\begin{thebibliography}{10}
\bibitem{Ciucu}
Mihai Ciucu, {\it Perfect matchings of cellular graphs}, J. Alg.
Combin. {\bf 5} (1996), 87--103.

\bibitem{EKLP}
 M. Elkies, G. Kuperberg, M. Larsen and J. Propp,
{\em Alternating sign matrices and domino tilings (Part~1)},
J. Alg\@. Combin\@. {\bf 1}
(1992), 111--132.
\bibitem{Propp}
J. Propp,
{\em The Many Faces of Alternating--Sign Matrices}, preprint.
%Submitted to Discrete Mathematics and Theoretical Computer Science.
\bibitem{shuffle}
J. Propp,
{\em Generalized Domino--Shuffling}, preprint.
\bibitem{yang}
B.--Y. Yang,
{\em Two Enumeration Problems about the Aztec Diamonds},
Ph.D. thesis, Mass. Inst. Tech., 1991.

\end{thebibliography}
\end{document}

%%% Local Variables: 
%%% mode: ams-tex
%%% TeX-master: t
%%% End: 